\documentclass[11pt]{amsart}
\usepackage{mathrsfs}
\usepackage{amssymb}
\usepackage{amsmath}
\usepackage{amsthm}

\makeatletter
\def\LaTeX{\leavevmode L\raise.42ex
    \hbox{\kern-.3em\size{\sf@size}{0pt}\selectfont A}\kern-.15em\TeX}
\makeatother

\newcommand{\BibTeX}{{\rm B\kern-.05em{\sc
          i\kern-.025emb}\kern-.08em\TeX}}

\makeatletter
\def\@currentlabel{2.1}\label{e:dispaa}
\def\@currentlabel{2.21}\label{e:dispau}
\def\@currentlabel{2.22}\label{e:dispav}
\def\@currentlabel{2.23}\label{e:dispaw}
\def\@currentlabel{2.24}\label{e:dispax}
\def\theequation{\thesection.\@arabic\c@equation}
\makeatother

\newcounter{mnotecount}[section]

\newcommand{\rmnote}[1]{}

\usepackage{amsfonts,latexsym}

\usepackage{graphics}

\renewcommand{\theequation}{\arabic{section}.\arabic{equation}}
\newcommand{\reset}{\setcounter{equation}{0}}

\newtheorem{thm}{Theorem}[section]

\newtheorem{prop}[thm]{Proposition}
\theoremstyle{definition}

\newtheorem{rem}[thm]{Remark}

\newcommand{\R}{\mathbb{R}}
\newcommand{\Sp}{\mathbb{S}}

\def\p{\partial}
\def \l{\lambda}
\def \e{\varepsilon}

\begin{document}
\title[weighted elliptic equation in exterior domains]{Existence and nonexistence results for a weighted elliptic equation in exterior domains}

\author{Zongming Guo}
\address{Department of Mathematics, Henan Normal University, Xinxiang, 453007, P.R. China}
\email{gzm@htu.cn}

\author{Xia Huang}
\address{Center for Partial Differential Equations, School of Mathematical Sciences, East China Normal
University, Shanghai 200241, China}
\email{xhuang@cpde.ecnu.edu.cn}

\author{Dong Ye}
\address{Center for Partial Differential Equations, School of Mathematical Sciences, East China Normal
University, Shanghai 200241, China}
\address{IECL, UMR 7502, D\'epartement de Math\'ematiques, Universit\'e de Lorraine, 57073 Metz, France}
 \email{dye@math.ecnu.edu.cn, dong.ye@univ-lorraine.fr}

\date{}
\begin{abstract}
We consider positive solutions to the weighted elliptic problem
\begin{equation*}
-\mbox{div} (|x|^\theta \nabla u)=|x|^\ell u^p \;\;\mbox{in $\R^N \backslash {\overline B}$},\quad
u=0 \;\; \mbox{on $\partial B$},
\end{equation*}
where $B$ is the standard unit ball of $\R^N$. We give a complete answer for the existence question for $N':=N+\theta>2$ and $p > 0$. In particular, for $N' > 2$ and $\tau:=\ell-\theta >-2$, it is shown that for $0< p \leq p_s:=\frac{N'+2+2 \tau}{N'-2}$, the only nonnegative solution to the problem is $u \equiv 0$. This nonexistence result is new, even for the classical case $\theta = \ell = 0$ and $\frac{N}{N-2} < p \leq \frac{N+2}{N-2}$, $N \geq 3$. The interesting feature here is that we do not require any behavior at infinity or any symmetry assumption.
\end{abstract}

\maketitle

\section{Introduction}
\setcounter{equation}{0}

We study positive solutions to the problem
\begin{equation}
\label{1.1}
-\mbox{div} (|x|^\theta \nabla u)=|x|^\ell u^p \;\;\mbox{in $\R^N \backslash {\overline B}$},\quad
u=0 \;\; \mbox{on $\partial B$}
\end{equation}
where $B=\{x \in \R^N, \; |x|<1\}$, $p > 0$.

\medskip
For $\theta=\ell=0$, the classical elliptic equation
\begin{equation}
\label{1.1c}
-\Delta u = u^p \;\;\mbox{in $\R^N \backslash {\overline B}$},\quad
u=0 \;\;\mbox{on $\partial B$}
\end{equation}
has been studied intensively, see for instance \cite{AMQ, AS, BM, DDMW, GS, Ni, Sa, Z} and the references therein.

\smallskip
It is well-known from \cite{AMQ, AS, GS, Ni} that \eqref{1.1c} does not admit any positive solution provided $0 < p \leq \frac{N}{N-2}$ for $N \geq 3$; or $p > 0$ for $N \leq 2$. It is also well-known from \cite{BM} that for any supercritical exponent $p>\frac{N+2}{N-2}$ and $N \geq 3$, the problem \eqref{1.1c} admits a unique positive radial solution $u$ satisfying
\begin{align*}
u(x) = O\left(|x|^{2-N}\right) \quad \mbox{as } |x| \to \infty.
\end{align*}
For $\frac{N}{N-2}<p \leq \frac{N+2}{N-2}$, $N \geq 3$, Theorem 2.2 in \cite{Sa} showed that \eqref{1.1c} does not admit any positive {\sl radial} solution. However, the following natural question remained open:
$$\mbox{Is there any {\sl non radial} positive solution of \eqref{1.1c} for } \frac{N}{N-2}<p \leq \frac{N+2}{N-2}, \; N \geq 3? \leqno{(*)}$$

It's worthy to mention two closely linked results, which seem to suggest that the answer for $(*)$ could be complex. Firstly, there exist infinitely many {\sl non radial} slow decay positive solutions to \eqref{1.1c} for any supercritical exponent $p$. More precisely, let $\Omega$ be a bounded smooth domain in $\R^N$ ($N \geq 3$), it was proved in \cite{DDMW} that for the homogenous Dirichlet problem $-\Delta u = u^p$ in $\R^N \backslash \Omega$ with $p>\frac{N+2}{N-2}$, there exist always infinitely many positive solutions such that
\begin{align*}
u(x) = O\left(|x|^{-\frac{2}{p-1}}\right) \quad \mbox{as } |x| \to \infty.
\end{align*}
Secondly, by Proposition 6.1 in \cite{Z}, Zhang showed that for any bounded Lipschitz domain $\Omega \subset \R^N$ ($N \geq 3$), $p > \frac{N+\ell}{N-2}$, the problem
\begin{equation*}
-\Delta u = |x|^\ell u^p \;\;\;\mbox{in $\R^N \backslash \Omega$},\quad u = f \;\;\;\mbox{on } \; \p\Omega
\end{equation*}
admits a positive solution, if $f$ is a {\it nontrivial} nonnegative function in $L^\infty(\p\Omega)$ and $\|f\|_\infty$ is small enough.

\medskip
In this paper, we study the existence of positive solutions to the more general equation \eqref{1.1}, under the following basic assumption:
\begin{equation}
\label{1.1-1}
N':=N+\theta>2,\;\;\;\tau:= \ell-\theta>-2.
\end{equation}
The elliptic problem like $-{\rm div}(a(x)\nabla u) = b(x)f(x, u)$ can be used to modelling some physical phenomena related to the equilibrium of continuous media, see \cite{DL}.

\medskip
We will show that
$$p_s:=\frac{N'+2+2 \tau}{N'-2}$$ is effectively the critical exponent for \eqref{1.1}. In particular, when $\theta=\ell = 0$, we prove that the answer to the question $(*)$ is indeed negative.

\medskip
Our main results are the follows.
\begin{thm}
\label{t1.2}
Let $N'>2$ and $\tau>-2$. The only nonnegative solution of \eqref{1.1} is $u_p \equiv 0$ provided
$0<p \leq p_s$.
\end{thm}

The interesting feature here is that we do not require any behavior at infinity or any symmetry assumption for $u$. Note that by scaling, our results hold obviously on $\R^N\backslash B_R$, $\forall \; R > 0$.

\smallskip
For $p \leq \frac{N'+\tau}{N'-2}$, we show a more general nonexistence result inspired by \cite{AMQ,AS}, see Theorem \ref{p3.2} below. For $1 < p \leq p_s$, we apply the moving-sphere method to get a monotonicity property for the eventual positive solution, then we conclude by contradiction with integral estimate and stability argument. Comparing to Zhang's result mentioned above, the homogenous Dirichlet boundary condition plays a crucial role for our nonexistence result whenever $\frac{N'+\tau}{N' - 2} < p \leq p_s$.

\begin{thm}
\label{t1.1}
Let $N'>2$, $\tau>-2$ and $p>p_s$. Then the problem \eqref{1.1} admits a unique positive radial solution $u_p \in C^2 (\R^N \backslash B)$.
\end{thm}

Combining with the approaches in \cite{AMQ, GW, St2, St}, we show finally a complete picture for the existence of positive solution to \eqref{1.1}, for all $N' > 2$ and $p > 0$.
\begin{thm}
\label{t4.1}
Let $N'>2$ and $p > 0$.
\begin{itemize}
\item[(1)] If $\tau > -2$, then \eqref{1.1} admits a positive solution if and only if $p > \frac{N'+2+2 \tau}{N'-2}$.
\item[(2)] If $\tau < -2$, then \eqref{1.1} admits a positive solution for any $p > 1$ and $p \in (0, 1)$.
\item[(3)] If $\tau < -2$ and $p = 1$, then \eqref{1.1} admits a positive solution if and only if $\lambda_1(L_1) = 0$ in $H^{1, \theta}_0(B)$ where
$$L_1(v) := -{\rm div}(|x|^\theta\nabla v) - |x|^{-(4+\tau-\theta)}v,$$
and $H^{1, \theta}_0(B)$ is the completion of $C^\infty_0(B)$ under the norm $\|u\| = \|\nabla u\|_{L^2(B, |x|^\theta dx)}$.
\item[(4)] Let $\tau = -2$. then \eqref{1.1} has a positive solution if and only if $p = 1$, $N' \geq 4$; or $p > 1$.
\end{itemize}
\end{thm}

In our analysis, we shall use the Kelvin transformation
\begin{align}
\label{Ke}
v(y)=|x|^{N'-2} u(x)\;\; \mbox{where } \; y=\frac{x}{|x|^2}.
\end{align}
If $u$ is a solution to \eqref{1.1}, then $v$ satisfies the problem
\begin{equation}
\label{1.2}
-\mbox{div} (|y|^\theta \nabla v)=|y|^{\sigma} v^p \;\;\mbox{in $B \backslash \{0\}$},\quad v=0 \;\; \mbox{on $\partial B$}.
\end{equation}
Here
$$\sigma=(N'-2) (p-1)-(4+\tau-\theta).$$
By direct calculation, there holds $\tau':=\sigma-\theta>-2$ if
$p>\frac{N'+\tau}{N'-2}$ and
\begin{equation}
\label{1.2-1}
p > p_s > 1\;\; \mbox{if and only if} \;\; 1 < p <p_s':=\frac{N'+2+2 \tau'}{N'-2}.
\end{equation}

Theorem \ref{t1.2} implies that problem \eqref{1.2} does not admit any nontrivial nonnegative solution in $C^2(B \backslash \{0\})$ for $N'>2$, $\tau'>-2$ and $p \geq p_s'$. A similar nonexistence result was proved in \cite{GW} with the additional assumption that $v \in C^0 ({\overline B})$. We prove Theorem \ref{t1.1} by using a positive radial solution $v_p \in C^2 (B \backslash \{0\}) \cap C^0 ({\overline B})$ to \eqref{1.2} given in \cite{GW}.

\section{Existence for supercritical case}
\setcounter{equation}{0}
Here we prove Theorem \ref{t1.1} by using the Kelvin transformation \eqref{Ke}, and we are looking for a solution $v$ satisfying the equation
\eqref{1.2}.

\medskip
Recall that $\sigma=(N'-2) (p-1)-(4+\tau-\theta)$ and $\tau'=\sigma-\theta>-2$. As $1<p<p_s'$ provided $p>p_s$, it follows from \cite{GW} that \eqref{1.2}
admits a positive radial solution $v_p \in C^2 (B \backslash \{0\})
\cap C^0 ({\overline B})$. Hence $r\mapsto v_p(r)$ satisfies the equation
\begin{equation}
\label{2.2}
-(r^{N'-1} v')'= r^{N'+\tau'-1} v^p(r) \;\; \mbox{in $(0,1)$}, \quad v(1) = 0.
\end{equation}
As $r^{N'-1} v_p'$ is decreasing in $r$, $N' > 2$; $v_p$ is positive and uniformly bounded in $(0, 1)$, $v_p$ must satisfy the following property:
\begin{align}
\label{v'0}
\lim_{r \to 0} r^{N'-1} v'(r)=0.
\end{align}
Then $v_p'(r)<0$ in $(0,1)$, this implies that $v_p(0)>0$.

\medskip
We now prove that $v_p$ is the unique positive radial solution of \eqref{1.2} in $C^2 (B \backslash \{0\}) \cap C^0 ({\overline B})$. Suppose that there is another such radial solution ${\underline v}$, let
$$\lambda^{\frac{2+\tau'}{p-1}}=\frac{{\underline v} (0)}{v_p(0)}, \quad \mbox{and }\; w(r):=\lambda^{\frac{2+\tau'}{p-1}} v_p(\lambda r) \;\; \mbox{for all } r \in [0,\lambda^{-1}].$$
Clearly, $w$ satisfies the equations
\begin{equation}
-(r^{N'-1} w')'= r^{N'+\tau'-1} w^p \;\; \mbox{in } (0,\lambda^{-1}), \quad w(0)= \underline v(0).
\end{equation}
Similarly as above,
\begin{equation*}
\lim_{r \to 0} r^{N'-1} w'(r)=0.
\end{equation*}

Consider $h = w - \underline v$, so $h(0) = 0$ and $h$ satisfies the following integral equation
\begin{equation}
\label{2.5}
h(r) = -\int_0^r t^{1-N'} \int_0^t s^{N'+\tau'-1}\Big[w^p(s) - {\underline v}^p(s)\Big] dsdt, \quad \forall\; 0 \leq r \leq \min\left(1, \lambda^{-1}\right).
\end{equation}
Fix any $0 < r_0 < \min(1, \lambda^{-1})$. By the continuity of $w$, $\underline v$ in $[0, r_0]$, the mean value theorem and Fubini's theorem imply that there exists $C > 0$ such that
$$|h(r)| \leq C\int_0^r s^{\tau'+1}|h(s)|ds, \quad \forall \; r \in [0, r_0].$$
Remark that $s^{\tau'+1} \in L^1_{loc}(0, \infty)$ since $\tau' > -2$, the classical Gronwall formula implies that $h \equiv 0$ in $[0, r_0]$. Therefore $w \equiv {\underline v}$ in $[0, r_0]$, hence in $[0, \min(1, \lambda^{-1}))$. This is possible if and only if $\lambda = 1$ seeing the Dirichlet boundary conditions and $v_p, {\underline v} > 0$ in $B$. To conclude, $\underline v \equiv v_p$ in $[0, 1]$.

\medskip
Conversely, let $u\in C^2 (\R^N \backslash B)$ be a positive radial solution to \eqref{1.1}, then
\begin{equation*}
-(r^{N'-1} u')'= r^{N'+\tau-1} u^p(r) \;\; \mbox{in $[1, \infty)$}.
\end{equation*}
So we have $r^{N'-1} u'$ is decreasing and $u(r) = O\big(r^{2 - N'}\big)$ at infinity. Under the Kelvin transform \eqref{Ke}, we get a bounded radial function $v \in C^2(\overline B\backslash\{0\})$ which satisfies \eqref{2.2}. As above, the boundedness of $v$ implies \eqref{v'0}, hence $v$ is decreasing by \eqref{2.2} and $\lim_{r \to 0^+} v(r)= \ell_0 \in \R$ exists, because $v$ is bounded. In other words, $v \in C^0(\overline B)$. The uniqueness of solution for \eqref{2.2} in $C^2 (B \backslash \{0\})
\cap C^0 ({\overline B})$ completes the proof. \qed

\section{Nonexistence for $1 < p \leq p_s$}
\setcounter{equation}{0}

Let $1 < p \leq p_s$, we show first a monotonicity property for positive solution of \eqref{1.1} if it exists.

\begin{prop}
\label{p3.1}
Let $N'>2$, $\tau>-2$, $1<p \leq p_s$ and $u \in C^2 (\R^N \backslash B)$ be a nontrivial nonnegative solution to \eqref{1.1}. Then
$|x|^{\frac{N'-2}{2}} u(x)$ is increasing with respect to the radius $r = |x|$.
\end{prop}

\noindent
{\bf Proof.} By the strong maximum principle, $u>0$ in $\R^N \backslash {\overline B}$. We introduce the following transform:
\begin{equation}
\label{3.1}
v(t, \omega)=r^{\frac{N'-2}{2}} u(r, \omega), \;\;\; t=\ln r.
\end{equation}
Then, $v$ satisfies the problem
\begin{equation}
\label{3.2}
\left \{ \begin{array}{l} v_{tt}+\Delta_{\Sp^{N-1}} v+e^{p_*t} v^p-\frac{(N'-2)^2}{4} v=0, \;\;\; (t, \omega) \in (0, \infty) \times \Sp^{N-1}, \\
v(0, \cdot)|_{\Sp^{N-1}}=0, \end{array} \right.
\end{equation}
where $$p_*=\frac{(N'+2+2 \tau)-(N'-2)p}{2} \geq 0.$$

Let $\Sigma:= (0, \infty) \times \Sp^{N-1}$. For $T > 0$, define
$$\Sigma_T:= (0, T) \times \Sp^{N-1}, \quad S_T= \{T\} \times \Sp^{N-1}.$$
Furthermore, we denote by $(t, \omega)_T$ the reflection of $(t, \omega)$ with respect to $S_T$, namely,
$$(t, \omega)_T =(2T-t, \omega).$$
Then the functions
$$v_T (t, \omega):=v((t, \omega)_T)=v(2T-t, \omega) \quad \mbox{and} \quad \xi_T=v-v_T$$
are well-defined in $\Sigma_T$ for $T>0$. Moreover, $\xi_T$ satisfies the following equation in $\Sigma_T$:
$$(\xi_T)_{tt}+\Delta_{\Sp^{N-1}} \xi_T+e^{p_*t} v^p-e^{p_*(2T-t)} v_T^p-\frac{(N'-2)^2}{4} \xi_T=0.$$
As $p_* \geq 0$, there holds $e^{p_* (2T-t)} \geq e^{p_*t}$ for any $t \leq T$. So we obtain
$$(\xi_T)_{tt}+\Delta_{\Sp^{N-1}} \xi_T+e^{p_* t} (v^p - v_T^p) -\frac{(N'-2)^2}{4} \xi_T \geq 0 \;\; \mbox{in $\Sigma_T$},$$
hence
\begin{equation}
\label{3.4}
(\xi_T)_{tt}+\Delta_{\Sp^{N-1}} \xi_T+\eta \xi_T \geq 0 \;\; \mbox{in $\Sigma_T$},\quad
\xi_T \leq 0 \;\;\mbox{on $\partial \Sigma_T$}
\end{equation}
where
$$\eta (t, \omega):=c(t, \omega)-\frac{(N'-2)^2}{4} \quad \mbox{and} \quad c(t, \omega) = pe^{p_*t}\psi^{p-1}(t, \omega),$$
with $\psi(t, \omega)$ in the interval formed by $v(t, \omega)$ and $v_T(t, \omega)$. The fact that $\xi_T \leq 0$ on $\partial \Sigma_T$ for any $T > 0$ comes easily from $u=0$ on $\partial B$.

\smallskip
By the continuity of $v$, we see that $\eta(t, \cdot)$ tends uniformly to $-\frac{(N'-2)^2}{4}$ on $\Sp^{N-1}$ as $T\to 0^+$. This implies that there exists $\delta_1 > 0$ such that
\begin{equation*}
\eta (t, \omega)<0 \;\; \mbox{for $(t, \omega) \in \Sigma_T$ and $0 < T \leq \delta_1$}.
\end{equation*}
Applying the maximum principle to \eqref{3.4}, it follows that $\xi_T \leq 0$ in $\Sigma_T$ provided $T \leq \delta_1$. As $\xi_T(0, \omega) < 0$, the strong maximum principle (see for example \cite{GT}) yields that for $0 < T \leq \delta_1$, there hold
\begin{equation}
\label{3.8}
\xi_T (t, \omega)<0 \;\; \mbox{in $\Sigma_T$}, \quad \frac{\partial \xi_T}{\partial t} (t, \omega)>0
\;\; \mbox{on $S_T$}.
\end{equation}
Remark also that
\begin{equation}
\label{3.7}
\frac{\partial \xi_T}{\partial t} (T, \omega) = 2 \frac{\partial v}{\partial t} (T, \omega), \quad \forall\; T > 0.
\end{equation}

Define now
\begin{equation*}
T_0=\sup \{T_1 > 0: \; \mbox{\eqref{3.8} holds for $T \leq T_1$} \}.
\end{equation*}
Then $T_0 \geq \delta_1$ is well-defined. We want to claim that $T_0 = \infty$.

\medskip
Suppose the contrary that $T_0<\infty$, then $\xi_{T_0} \leq 0$ in $\Sigma_{T_0}$ by continuity. Seeing the elliptic equation \eqref{3.4}, as $\xi_{T_0}(0, \omega) < 0$, the strong maximum principle means that \eqref{3.8} holds true for $T = T_0$. Using \eqref{3.7}, we deduce that $\p_t v > 0$ in $(0, T_0]\times \Sp^{N-1}$. Recall that $u \in C^1(\R^N\backslash B)$, by the compactness of $\Sp^{N-1}$, there exists $\e_1 \in (0, T_0)$ such that
$$\frac{\partial v}{\partial t} (t, \omega) > 0 \quad \mbox{for }(t, \omega) \in (0, T_0+ 3\e_1] \times \Sp^{N-1}.$$
We get then
\begin{align}
\label{new2}
\xi_T(t, \omega) < 0 \;\; \mbox{in } \Sigma_T\backslash \Sigma_{T_0 - \e_1}, \quad \forall\; T_0 \leq T \leq T_0 + \e_1.
\end{align}
On the other hand, $\max_{\overline \Sigma_{T_0 - \e_1}} \xi_{T_0} < 0$ by continuity and \eqref{3.8} with $T = T_0$. Using the continuity of $v$, there exists $\e_2 > 0$ such that
$$v(t, \omega) - v(2T - t, \omega) < 0\;\; \mbox{in } \overline \Sigma_{T_0 - \e_1}, \quad \forall\; T_0 \leq T \leq T_0 + \e_2.$$
Combining with \eqref{new2}, let $\e = \min(\e_1, \e_2) > 0$, there holds $\xi_T \leq 0$ in $\Sigma_T$ for any $T \in [T_0, T_0 + \e]$.

\medskip
Again, as $\xi_T < 0$ on $\{0\}\times \Sp^{N-1}$, the strong maximum principle applying to \eqref{3.4} shows that \eqref{3.8} remains valid for $T \in [T_0, T_0 + \e]$, this is a contradiction to the definition of $T_0$. Therefore $T_0$ must be $\infty$, hence \eqref{3.8} holds true for any $T > 0$.

\smallskip
Using \eqref{3.7}, for any $\omega \in \Sp^{N-1}$, the function $r\mapsto r^{\frac{N'-2}{2}}u(r,\omega)$ is increasing in $[0, \infty)$.
\qed

\medskip
Now we are ready to prove Theorem \ref{t1.2} for $1 < p \leq p_s$. Suppose that \eqref{1.1} admits a nontrivial nonnegative solution $u \in C^2 (\R^N \backslash B)$. We claim that
$$\varphi(x) :=\nabla u \cdot x+\frac{2+\tau}{p-1} u(x)>0  \;\; \mbox{in }\; \R^N \backslash \overline B.$$
In fact, as $1< p \leq p_s$, we have
$$
|x|^{\frac{2+\tau}{p-1}} u(x)=|x|^{\beta}|x|^{\frac{N'-2}{2}}u(x),\quad \text{where }\; \beta:=\frac{2+\tau}{p-1}-\frac{N'-2}{2}\geq 0.
$$
So the function $|x|^{\frac{2+\tau}{p-1}}u(x)$ is increasing in $|x|$, thanks to Proposition \ref{p3.1}. By the above proof, $\p_t v > 0$ in $(0, \infty) \times \Sp^{N-1}$, hence $\varphi > 0$ in $\R^N \backslash \overline B$.

\medskip
Direct calculation via scaling yields that $\varphi \in C^1(\R^N \backslash B)$ is a positive weak solution to the linearized equation
\begin{equation*}
L_u(w) := -\mbox{div}(|x|^\theta \nabla w) - p |x|^\ell u^{p-1} w = 0 \;\;\; \mbox{in $\R^N \backslash {\overline B}$}.
\end{equation*}
By classical elliptic theory, for any bounded smooth domain $\Omega \subset \R^N \backslash {\overline B}$,
\begin{equation}
\label{3.16}
\lambda_{1,\Omega} (L_u)>0,
\end{equation}
where $\lambda_{1,\Omega} (L_u)$ stands for the first eigenvalue of the Dirichlet problem
$$L_u(h) =\lambda h \;\; \mbox{in }\; \Omega, \quad
h=0 \;\; \mbox{on }\; \partial \Omega.$$
So $u$ is stable in $\R^N \backslash {\overline B}$.

\medskip
Using Theorem 3.1 in \cite{DG}, we have
\begin{equation}
\label{3.17}
\int_{\R^N \backslash {\overline B}} \Big(|x|^\theta |\nabla u|^2+|x|^\ell u^{p+1} \Big)dx<\infty.
\end{equation}
However, Proposition \ref{p3.1} implies then the existence of $C > 0$ such that
$$u(x) \geq C |x|^{-\frac{N'-2}{2}} ,\;\; \forall\;|x| \geq 2.$$
Recall that $1< p\leq p_s$, there holds
\begin{equation*}
\int_{|x|>2}|x|^\ell u^{p+1} dx \geq C \int_2^\infty r^{N-1+\ell} r^{-\frac{(p+1)(N'-2)}{2}} dr=\infty,
\end{equation*}
since
$$N-1 + \ell -\frac{(p+1)(N'-2)}{2} \geq N'+ \tau - 1 - \frac{(p_s+1)(N'-2)}{2} = -1.$$
This contradiction with \eqref{3.17} means that there is no nontrivial nonnegative solution to \eqref{1.1}. \qed

\section{Further results}
\reset
First, we consider the existence and nonexistence of positive solution to
\begin{equation}
\label{3.19}
-\mbox{div} (|x|^\theta \nabla u) \geq |x|^\ell f(u)\;\;\; \mbox{in $\R^N \backslash {\overline B}_R$, $R > 0$.}
\end{equation}
Remark that we do not impose any boundary condition or the value of $R$.

\smallskip
A first nonexistence result has been obtained by Gidas-Spruck (see Theorem A.3 in \cite{GS}), they showed that for $\theta=\ell$, $N'= N+\theta >2$ and $f(u)=u^p$, \eqref{3.19} admits only trivial non-negative solution $u\equiv 0$ provided $1<p\leq \frac{N'}{N'-2}$. Since then, many existence and nonexistence results have been established for \eqref{3.19} under suitable conditions, interested readers can look at \cite{AMQ, AS, LSS, St2, St, Z} and the references therein.

\medskip
Inspired by Theorem 8 of \cite{AMQ}, we show a very general result to characterize the existence of positive solution to \eqref{3.19}. In particular, the following result implies Theorem \ref{t1.2} for $0 < p \leq \frac{N'+\tau}{N' - 2}$.
\begin{thm}
\label{p3.2}
Assume that $N'>2$, $\tau>-2$ and $f: (0, \infty) \to (0, \infty)$ is continuous.
Then the problem \eqref{3.19} admits a positive classical solution for some $R > 0$ if and only if
\begin{equation}
\label{3.20}
\int_0^\delta f(t)t^{-\frac{2(N'-1)+\tau}{N'-2}} dt < \infty,
\end{equation}
for some $\delta>0$.
\end{thm}

\noindent
{\sl Proof.} We need only to adapt the approach in \cite{AMQ}. Note that to prove Lemma 4 in \cite{AMQ}, the key points are the maximum principle and the rotational invariance of the equation, which obviously remain true in our setting.

\smallskip
Therefore, we can claim that if \eqref{3.19} admits a positive solution, then for some $R_0> R$, there exists a $C^1$ positive radial function ${\overline u}$ verifying
$$-\mbox{div} (|x|^\theta \nabla u)=|x|^\ell f(u) \;\;\; \mbox{in $\R^N \backslash B_{R_0}$}$$
in the weak sense. Moreover, if ${\overline u}(x)=v(|x|)$, then $v \in C^2 (R_0, \infty)$ and there exists
$R_1 \geq R_0$ such that $v$ is monotone for $r>R_1$. Note that $v$ satisfies
\begin{equation}
\label{3.21}
-v''-\frac{N'-1}{r} v'=r^\tau f(v) \;\; \mbox{in } \; [R_1, \infty).
\end{equation}
As $N' > 2$, there holds $ v(r)\leq -C_1r^{2-N'} + C_2$. Hence $v$ is bounded in $[R_1, \infty)$ and $\lim_{r\to\infty} v(r) = \kappa\geq 0$ exists.

\medskip
If $\kappa > 0$, there exists $C> 0$ such that
 $(r^{N'-1}v')' \leq -Cr^{N'+\tau -1}$ for $r$ large enough. Using $N'+\tau > 0$ and integrating, we deduce that $v'(r)\leq -Cr^{\tau + 1}$ for large $r$. As $\tau > -2$, this implies that $\lim_{r\to\infty}v(r)
 = -\infty$ which is impossible. Hence we can claim that $v'(r)<0$ in $[R_1, \infty)$ and $\lim_{r \to \infty} v(r)=0$.

\medskip
Similarly to \cite{AMQ}, we introduce the change of variables
\begin{align}
\label{w}
s := r^{2 - N'}, \quad w(s)= v(r).
\end{align}
Then $w$ is a positive solution of
\begin{equation}
\label{3.22}
-w''(s)=a s^{-\gamma} f(w) \;\;\mbox{in $(0, s_0)$},\quad w(0)=0,
\end{equation}
for some positive constants $a$ and $s_0$, where (as $\tau > -2$)
$$\gamma=\frac{2(N'-1)+\tau}{N'-2}>2.$$
Finally, applying Theorem 6 in \cite{AMQ}, a positive solution of \eqref{3.22} exists if and only if
\eqref{3.20} holds true.\qed

\begin{rem}
Let $N'> 2$, $\tau>-2$, then
$$p-\frac{2(N'-1)+\tau}{N'-2} \leq -1 \quad \Leftrightarrow \quad p \leq \frac{N'+\tau}{N'-2}.$$ For $f(t) = t^p$ with $0< p \leq \frac{N'+\tau}{N'-2}$, Proposition \ref{p3.2} implies that \eqref{3.19} cannot have positive solution for any $R > 0$, hence the equation \eqref{1.1} does not admit positive solution.
\end{rem}

\begin{rem}
It is well known that if $N' > 2$, $\tau > -2$ and $p > \frac{N'+\tau}{N'-2}$, $v(x) = Cr^{-\frac{\tau + 2}{p - 1}}$ with suitable $C(N', \tau, p) > 0$ will satisfy $-\mbox{div} (|x|^\theta \nabla v) = |x|^\ell v^p$ in $\R^N\backslash\{0\}$.
\end{rem}

Moreover, Przeradzki and Sta\'nczy \cite{St2, St} obtained that if $N > 2$, $\theta = 0$ and $\ell < -2$, the equation \eqref{1.1} admits a positive solution for any $p > 1$ and $p \in (0, 1)$, see Remark 5 in \cite{St} and Corollary 1 in \cite{St2}. It's not difficult to be convinced that their approach via radial solutions remains valid to \eqref{1.1} for any $N' > 2$ and $\tau < -2$. We get then a complete answer to the existence question of positive solution of \eqref{1.1}, under the only assumption $N' := N + \theta > 2$.

\medskip
\noindent
{\bf Proof of Theorem \ref{t4.1}.} The point (1) is given by Theorems \ref{t1.2}-\ref{t1.1}. The point (2) is ensured by \cite{St2, St} as explained above.

\medskip
Consider $\tau < -2$ and $p = 1$, assume that a positive solution exists to \eqref{1.1}. As the equation is linear, the existence of a positive solution to \eqref{1.1} yields a positive {\sl radial} solution $u$ to \eqref{1.1}, and there holds $u(r) = O(r^{2-N'})$ at infinity. By the Kelvin transform \eqref{Ke}, we get a solution $v$ to \eqref{2.2} with $p =1$. As above, there holds $v \in C^0(\overline B)$ and $v$ satisfies \eqref{v'0}. Using the equation, we see that $v \in H^{1, \theta}_0(B)$ is a positive solution to $L_1(v) = 0$, with $\tau' = -(4+\tau) > -2$. It means that $0$ is the principal eigenvalue to $L_1$ in $H^{1, \theta}_0(B)$.

\medskip
Fix now $\tau = -2$.
\begin{itemize}
\item If $0< p < 1$, we will prove that the problem \eqref{3.19} has no solution for any $R > 0$. Suppose the contrary, as in \cite{AMQ}, we get a positive {\sl radial} solution to ${\rm -div}(|x|^\theta\nabla v) = |x|^{\ell}v^p$ in $\R^N\backslash B_{R_1}$ for some $R_1 > R$. Under the change of variable \eqref{w}, there exists $a, s_0 > 0$ such that
\begin{align*}
-w'' = as^{-2}w^p \;\; \mbox{in } (0, s_0), \quad w(0) = 0.
\end{align*}
However, this is impossible seeing the proof of Theorem 8 (b) in \cite{AMQ}, where we can replace $N$ just by $N' > 2$.
\item If $p = 1$, as before, the existence of a positive solution to \eqref{1.1} is equivalent to the existence of a positive {\sl radial} solution. Assume that $u$ is a radial solution to \eqref{1.1}, let $w(s) = u(e^s)$, then $w$ is a solution to
$$w'' + (N' - 2)w' + w = 0 \;\; \mbox{in } (0, \infty), \quad w(0) = 0, \;\; \lim_{s\to\infty} w(s) = 0.$$
Clearly, such a positive solution $w$ exists if and only if the equation $\l^2 + (N'-2)\l + 1 = 0$ has only negative roots, or equivalently, $N' \geq 4$.
\item Finally let $p > 1$. Under the Kelvin transform, we are led to consider \eqref{1.2}. We check readily that $\tau' = \sigma - \theta = (N' - 2)(p-1)-2 > -2$ and $p$ is subcritical, since here
\begin{align*}
p < \frac{N' + 2 + 2\tau'}{N'-2} \quad \Leftrightarrow \quad (p-1)(N' - 2) > 0
\end{align*}
There exists a positive radial solution of \eqref{1.2} as in section 2, using the result in \cite{GW}. \hfill \qed
\end{itemize}

\begin{rem}
\label{rem2}
For $N'=2$, $\tau>-2$, similarly to \cite{AMQ}, we can show that the problem \eqref{3.19} admits a positive solution if and only if there exists $a>0$ such that $e^{at} f(t)$ is integrable at $+\infty$. Consequently, in this case, for any $p \in \R$, the equation \eqref{3.19} does not admit positive solution.
\end{rem}

\bigskip
\noindent
{\small {\bf Acknowledgements.} Z.G is supported by NSFC (No.~11571093 and 11171092). X.H is supported by NSFC (No.~11701181). D.Y is partially supported by Science and Technology Commission of Shanghai Municipality (STCSM) (No.~18dz2271000). The authors thank the anonymous referees for their careful reading.

\end{document}